# About one method of construction of interpolation trigonometric splines

Volodymyr Denysiuk, Dr of Phys-Math. sciences, Professor, Kiev, Ukraine

National Aviation University

kvomden@nau.edu.ua

## Annotation

The method of constructing spline classes in the form of trigonometric Fourier series whose coefficients have a certain decreasing order are considered. in turn, this decrement determines the number of continuous derivatives of sum of this series. By grouping members of this series according with the effect of overlaying and introducing a multiplier that provides the interpolation properties of the sum of these series on even grids, we obtain classes of trigonometric interpolation splines. Depending on the types of convergence factors with a certain decreasing order, different classes of such splines are obtained. The classes of trigonometric splines include classes of periodic polynomial splines of even and odd power; At the same time, there exist trigonometric splines that do not have polynomial analogues. An example of the construction of trigonometric splines is given.

## Introduction

As known, interpolation splines in the general case are called interpolation functions that have certain properties of smoothness. The most known at present are splines, which are sewn from polynomial polynomials. Appearing in the middle of the twentieth century, polynomial splines immediately attracted the attention of specialists in numerical methods due to their approximation properties. Thus, in particular, it was shown that polynomial splines are the best linear apparatus for the approximation of classes of smooth functions.

A certain generalization of polynomial splines is trigonometric splines proposed by J. Shoenberg's, which are sewn from the pieces of trigonometric polynomials. The theory of trigonometric splines is quite developed and has a number of branches.

However, polynomial and trigonometric splines (and their generalizations), which are composite functions, have a number of shortcomings that will be discussed below. Such situation makes the actual task of developing other approximation devices, which would have the properties of polynomial splines, and at the same time serve as the only analytical expression throughout the interval of the approximation of function. For developing one of these approximation devices, we used the Fourier trigonometric series, for which, as is known, exists a direct relationship between the differential properties of the sum of such series with the order of decreasing of coefficients of these series. To provide the same interpolation properties of the sum of such series in the role of their coefficients, the products of the coefficients of interpolation trigonometric polynomials with different convergence factors with certain decreasing orders were used. Some results of research in this direction are presented in this paper.

## Literature review

One of the first monographs on the theory of polynomial splines is a work [1]. Various aspects of the theory of such splines were considered, for example, in papers [2], [3], [4]. The approximative properties of splines are devoted to work [5].

Splines sewn from trigonometric polynomials were first described in [6]. One-dimensional trigonometric splines are bits-trigonometric polynomials, which are given between nodes of interpolation grids; they are a natural generalization of polynomial splines. Some theoretical questions related to the construction of trigonometric splines are covered, for example, in [7]. Together with splines wich are sewn from trigonometric curves, splines that are cross-linked with other curves of the standard type, such as the Overhauser splines and Bezier splines, are often considered. [8].

Different variants of construction of trigonometric splines of one and two variables are considered in [9], [10], [11], [12].

Some applications of trigonometric splines are considered in [13] [14].

As we have already said, both polynomial and trigonometric splines have a number of disadvantages. The most important of these is that these splines have a composite structure. This shortcoming significantly restricts the use of polynomial splines in many numerical analysis tasks; In addition, due to the composite structure, in practice, predominantly polynomial splines of the third degree and trigonometric splines, which sew a small number of derivatives, are used in practice.

## The purpose of the work

The purpose of the work is to develop a constructive method for constructing trigonometric splines that have a certain number of continuous derivatives and interpolate the function in nodes of uniform grids on interval $[0, 2\pi)$; this method is based on the properties of trigonometric Fourier series, for which exists a direct connection between the differential properties of the sums of such series with the order of decreasing of the coefficients of these series.

## Main Part

On the gap $[0, 2\pi)$ given a continuous function $f(t)$. Let also be given even grids on this gap $\Delta_N^{(0)} = \{t_i^{(0)}\}_{i=1}^N$, $t_i^{(0)} = \frac{2\pi}{N}(i-1)$, and $\Delta_N^{(1)} = \{t_i^{(1)}\}_{i=1}^N$, $t_i^{(1)} = \frac{\pi}{N}(2i-1)$  $N = 2n+1$, $n = 1, 2, \ldots$ . To shorten the recording, we introduce the grid indicator $I = 0, 1$; then both grids can be labeled as $\Delta_N^{(I)}$. Let's denote as $\{f(t_i^{(I)})\}_{i=1}^N = \{f_i^{(I)}\}_{i=1}^N$ set of function values $f(t)$ in the nodes of the grid $\Delta_N^{(I)}$. Consider a trigonometric polynomial

$$T_n^{(I)}(t) = \frac{a_0^{(I)}}{2} + \sum_{k=1}^n a_k^{(I)} \cos kt + b_k^{(I)} \sin kt, \quad (1)$$

which interpolates the function $f(t)$ on the grid $\Delta_N^{(I)}$. It is known (see [17]) that the coefficients of this polynomial are determined by the formulas

$$a_k^{(I)} = \frac{2}{N} \sum_{j=1}^N f_j^{(I)} \cos kt_j^{(I)}, \quad b_k^{(I)} = \frac{2}{N} \sum_{j=1}^N f_j^{(I)} \sin kt_j^{(I)},$$
$$k = 0, 1, \ldots, n; \quad\quad\quad\quad k = 1, 2, \ldots, n. \quad (2)$$

We put in accordance the trigonometric polynomial $T_n^{(I)}(t)$ a generalized trigonometric function [15]

$$ST^{(I)}(\nu, r, N, t) = \frac{a_0^{(I)}}{2} + \sum_{k=1}^\infty \nu_k(r) \left[ a_k^{(I)} \cos kt + b_k^{(I)} \sin kt \right], \quad (3)$$

Where multipliers $\nu = \{\nu_k(r)\}_{k=1}^\infty$ have a descending order $O(k^{-(1+r)})$, $r > 0$, and cooficients $a_k^{(I)}, b_k^{(I)}$ are calculated by the formulas (2) for any one values $k$ ($k = 1, 2, \ldots$).

It is clear that the series (3) constructed in this way, on the basis of the Weierstras's sign, coincides evenly and its sum is a periodic function, continuous on the entire numerical axis. In the future, wanting to deal with ordinary derivatives, we will consider the case of only natural values of the parameter, that is, we will assume that this parameter takes value $r = 1, 2, \ldots$ . It is easy to see that this is a periodic function $ST^{(I)}(\nu, r, N, t)$ is continuous and has a continuous derivative order $r - 1$, that is $ST^{(I)}(\nu, r, N, t) \in C_{(-\infty, \infty)}^{r-1}$.

Now, we will require a function $ST^{(I)}(\nu,r,N,t)$ interpolated trigonometric polynomial (1) (and therefore and function $f(t)$) in the grid's nodes $\Delta_N^{(I)}$. To do this, we introduce the expression (3) multiplier $H_k^{(I)}(r,N)^{-1}$ which will be chosen from the conditions of interpolation.

Calculating the value of the function $ST^{(I)}(\nu,r,N,t)$ in the grid's nodes $\Delta_N^{(I)}$, due to the known effect of overlay, we get

$$ST^{(I)}(\nu,r,N,t_i^{(I)}) = \frac{a_0^{(I)}}{2} +$$
$$+ \sum_{k=1}^{n} \frac{1}{H_k^{(I)}(r,N)} \left[ a_k^{(I)} C_k^{(I)}(\nu,r,N,t_i^{(I)}) + b_k^{(I)} S_k^{(I)}(\nu,r,N,t_i^{(I)}) \right]. \quad (4)$$

where

$$C_k^{(I)}(\nu,r,N,t) = \left\{ \nu_k(r)\cos kt + \sum_{m=1}^{\infty} \left[ \nu_{mN+k}(r)\cos(mN+k)t + \nu_{mN-k}(r)\cos(mN-k)t \right] \right\},$$

$$S_k^{(I)}(\nu,r,N,t) = \left\{ \nu_k(r)\sin kt + \sum_{m=1}^{\infty} \left[ \nu_{mN+k}(r)\sin(mN+k)t - \nu_{mN-k}(r)\sin(mN-k)t \right] \right\}.$$

There is a theorem.

**Theorem 1.** Function $ST^{(I)}(\nu,r,N,t)$ interpolates a trigonometric polynomial $T_n^{(I)}(t)$ in the grid's nodes $\Delta_N^{(I)}$ only when multiplier $H_k^{(I)}(r,N)$ is determined by the formulas

$$H_k^{(I)}(r,N) = \left\{ \nu_k(r) + \sum_{m=1}^{\infty} (-1)^{mI} \left[ \nu_{mN+k}(r) + \nu_{mN-k}(r) \right] \right\}, \quad (5)$$

The proof of theorem is easy to obtain by comparing the expression (4) with (1) in the grid's nodes $\Delta_N^{(I)}$.

Given (5), the function $ST^{(I)}(\nu,r,N,t)$, can be written as

$$ST^{(I)}(\nu,r,N,t) = \frac{a_0^{(I)}}{2} + \sum_{k=1}^{n} H_k^{(I)}(r)^{-1} \left[ a_k^{(I)} C_k(\nu,r,N,t) + b_k^{(I)} S_k(\nu,r,N,t) \right]. \quad (7)$$

So, we built the function $ST^{(I)}(\nu,r,N,t)$, which, firstly, have certain differential properties, and secondly, interpolate the given function in the grid's nodes $\Delta_N^{(1)}$. These properties of functions $ST^{(I)}(\nu,r,N,t)$ give reason to call them **trigonometric splines of the order** $r$.

Let's consider some properties of trigonometric splines.

**Theorem 2.** Equality holds

$$\lim_{r \to \infty} ST^{(I)}(\nu,r,N,t) = T_n^{(I)}(t)$$

**Theorem 3.** Equality holds

$$\lim_{N \to \infty} ST^{(I)}(\nu,r,N,t) = T_n^{(I)}(t)$$

Proof we will carry out for a function $ST^{(0)}(\nu,r,N,t)$. Give (8) in the form

$$ST^{(0)}(\nu,r,N,t) = \frac{a_0^*}{2} + \sum_{k=1}^{n} \left[ a_k^* \frac{C_k(\nu,r,N,t)}{H_k^{(0)}(r)} + b_k^* \frac{S_k(\nu,r,N,t)}{H_k^{(0)}(r)} \right].$$

Let's show that

$$\lim_{r \to \infty} \frac{C_k(\nu,r,N,t)}{H_k^{(0)}(r)} = \cos kt \; ; \quad \lim_{r \to \infty} \frac{S_k(\nu,r,N,t)}{H_k^{(0)}(r)} = \sin kt \; .$$

Really because

$$C_k(\nu,r,N,t) = \nu_k(r)\cos kt + \sum_{m=1}^{\infty} \nu_{mN+k}(r)\cos(mN+k)t + \nu_{mN-k}(r)\cos(mN-k)t \; ;$$

$$S_k(\nu,r,N,t) = \nu_k(r)\sin kt + \sum_{m=1}^{\infty} \nu_{mN+k}(r)\sin(mN+k)t - \nu_{mN-k}(r)\sin(mN-k)t \; ,$$

$$H_k^{(0)}(r,N) = \left\{ v_k(r) + \sum_{m=1}^{\infty} (-1)^{mI} \left[ v_{mN+k}(r) + v_{mN-k}(r) \right] \right\},$$

then, taking into account the limitation of trigonometric functions, as well as the fact that the parameter $r$ determines the order of decreasing values $v_k(r)$ with the growth of their index, we arrive at the conclusion that at the main part of the expressions $C_k(v,r,N,t)$, $S_k(v,r,N,t)$ та $H_k(r,N)$ are the first members of the sums that determine these expressions. From this follows the statement of the theorem.

The second theorem proves similarly.

As we said above, the multipliers $v = \{v_k(r)\}_{k=1}^{\infty}$ have a decreasing order $O(k^{-(1+r)})$. In this paper, we restrict ourselves to considering the case when these factors have the form

$$v1_k(r) = \left[\operatorname{sinc}\left(\frac{\pi}{N}k\right)\right]^{1+r}; \quad v2_k(r) = \left|\operatorname{sinc}\left(\frac{\pi}{N}k\right)\right|^{1+r}; \quad v3_k(r) = \left(\frac{1}{k}\right)^{1+r},$$

where $\operatorname{sinc}(x) = \sin(x)/x$.

First of all it's clear that the factors $v1_k(r)$ і $v2_k(r)$ coincide at odd values of the parameter $r$, and, consequently lead to the same results. Next, direct calculations are easy to verify that the factors $v2_k(r)$ і $v3_k(r)$ lead to identical results for all natural values of the parameter $r$.

Thus, at odd values of the parameter, all three types of factors lead to identical results.

For even values of the parameter $r$ the results obtained with the multipliers $v1_k(r)$ and $v2_k(r)$ are different; but with increasing parameters $r$ (or $N$) these differences are directed to $0$, according to Theorems 2, 3.

An interesting question is the connection between the classes of trigonometric splines constructed in this way with the classes of periodic simple polynomial splines.

In [17] it is shown that periodic simple polynomial splines of odd power $2k-1$ ($k=1,2,...$) coincide with trigonometric splines of order $r = 2k-1$ for all types of multipliers on a grid $\Delta_N^{(0)}$.

Periodic simple polynomial splines of the even degree $2k$ ($k=1,2,...$) coincide with trigonometric splines of the same order in a grid $\Delta_N^{(1)}$ for all types of factors.

This fact allows us to extend to the listed trigonometric splines all the results of estimates of the approximation of functions obtained for polynomial splines as even and odd power. For example, it can be states that such trigonometric splines are the best linear apparatus for the approach of classes of smooth functions [5].

It is easy to come to this conclusion in another way, using the results of works [15], [16], in which, by means of factors $v1_k(r)$ polynomial $B$ - splines were built..

Note that polynomial analogs of trigonometric splines of even order on a grid $\Delta_N^{(0)}$ for all types of factors $v1_k(r)$, $v2_k(r)$ і $v3_k(r)$ are unknown for author; Similarly, unknown polynomial analogues and trigonometric splines of odd order on a grid $\Delta_N^{(1)}$ for values $r = 2r+1$, ($k=1,2,...$) for all types of factors $v1_k(r)$, $v2_k(r)$ і $v3_k(r)$. For the case of $r = 1$ polynomial analog of trigonometric spline of odd order on a grid $\Delta_N^{(1)}$ was built in[17].

It follows from this that the class of trigonometric interpolation splines is wider than the class of simple polynomial splines; In other words, the class of simple polynomial splines is a subset of a class of trigonometric splines.

Note that a certain disadvantage of trigonometric splines is their frequency. However, taking into account the known methods of phantom functions and phantom nodes of the periodic continuation of nonperiodic functions [16], we can conclude that periodicity does not play a decisive role when comparing polynomial and trigonometric splines.

Consider a numerical example.

Let's put $N = 9$ and set the interval $[0, 2\pi)$ of grid $\Delta_9^{(I)} = \{t_i^{(I)}\}_{i=1}^9$, $t_i^{(0)} = \frac{2\pi}{9}(i-1)$, and $t_i^{(1)} = \frac{\pi}{9}(2i-1)$. We also set the nodes of the grid $\Delta_N^{(I)}$ the value of function $\{f(t_i^{(I)})\}_{i=1}^9 = \{2,1,3,2,4,1,3,1,3\}$.

Let's put it for certainty $I = 0$. Trigonometric polynomial $T_n^{(0)}(t)$, which interpolates these function values in the nodes of the grid $\Delta_9^{(0)}$, shown in fig 1.

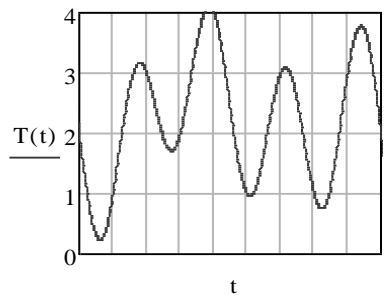

*Fig.1.* Graph of interpolation trigonometric polynomial $T_n^{(0)}(t)$

Graphs of trigonometric spline $ST^{(0)}(v1, r, 9, t)$ at different values of the parameters $r$ are shown in fig. 2.

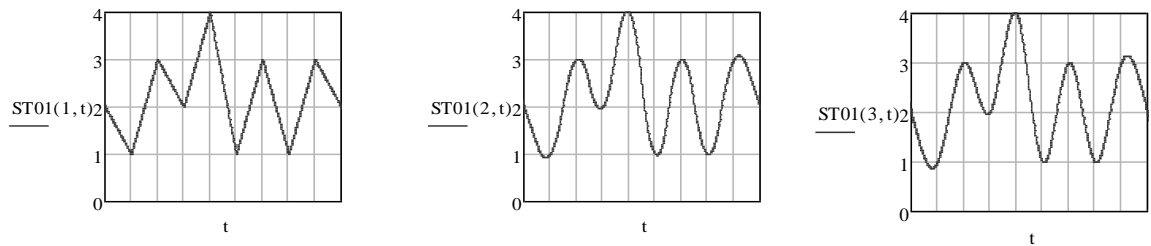

*Fig. 2.* Trigonometric spline $ST^{(0)}(v1, r, 9, t)$ at parameter values $r = 1, 2, 3$.

As we already said, with odd values of the parameter $r$ splines $ST^{(0)}(v1, r, 9, t)$ $ST^{(0)}(v2, r, 9, t)$ and $ST^{(0)}(v3, r, 9, t)$ coincides; so we do not provide graphs for these splines.

For even values of the parameter $r$ splines $ST^{(0)}(v2, r, 9, t)$ and $ST^{(0)}(v3, r, 9, t)$ coincide; but a spline $ST^{(0)}(v1, r, 9, t)$ are differs from them; Considering this, we give only the graphs of deviation $E0(r, t) = ST^{(0)}(v1, r, 9, t) - ST^{(0)}(v3, r, 9, t)$ for some even values of this parameter.

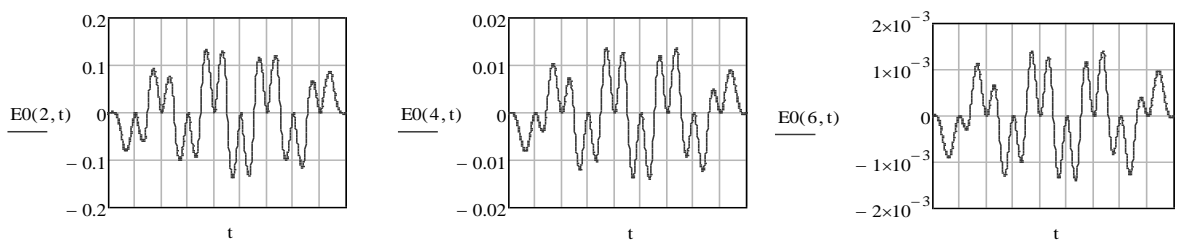

*Fig3.* Deviation of spline $ST^{(0)}(v1, r, 9, t)$ from spline $ST^{(0)}(v3, r, 9, t)$

with parameter values $r = 2, 4, 6$.

Now we consider the case $I = 1$. Trigonometric polynomial $T_n^{(1)}(t)$, which interpolates these function values in the nodes of the grid $\Delta_9^{(1)}$, is shown in fig 4.

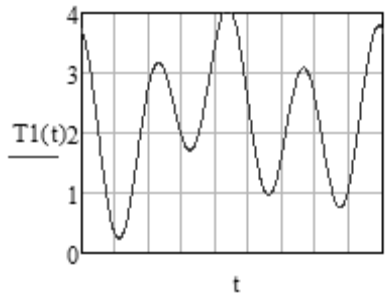

*Fig.4. Graph of interpolation trigonometric polynomial $T_n^{(1)}(t)$*

Graphs of trigonometric spline $ST^{(1)}(v1, r, 9, t)$ at different values of parameters $r$ are shown in Fig. 5.

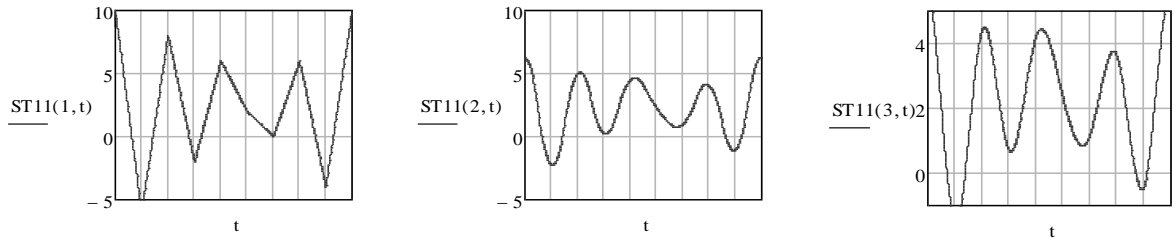

*Fig.. 5. Trigonometric spline $ST^{(1)}(v1, r, 9, t)$ with parameter values $r = 1, 2, 3$.*

As before, with odd values of the parameter $r$ splines $ST^{(1)}(v1, r, 9, t)$, $ST^{(1)}(v2, r, 9, t)$ and $ST^{(1)}(v3, r, 9, t)$ coincide; so spline graphs $ST^{(1)}(v2, r, 9, t)$ and $ST^{(1)}(v3, r, 9, t)$ we do not provide.

For even values of parameter $r$ splines $ST^{(1)}(v2, r, 9, t)$ and $ST^{(1)}(v3, r, 9, t)$ coincide; but a spline $ST^{(1)}(v1, r, 9, t)$ differs from them; Considering this, we give only the graph of deviation $E1(r, t) = ST^{(1)}(v1, r, 9, t) - ST^{(1)}(v3, r, 9, t)$ for some even values of this parameter.

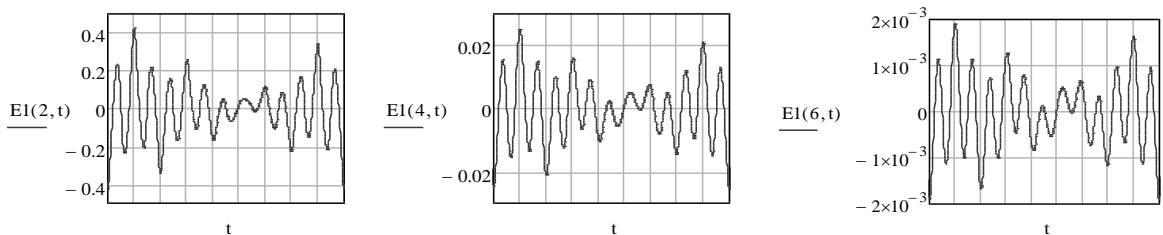

*Fig.6. Deviation of spline $ST^{(1)}(v1, r, 9, t)$ from spline $ST^{(1)}(v3, r, 9, t)$ While values are $r = 2, 4, 6$.*

We end with consideration of the example.

# Conclusions

1. A constructive method for constructing trigonometric splines with a certain number of continuous derivatives and interpolating the function in nodes of uniform grids in the space $[0, 2\pi)$ is developed; this method is based on the properties of trigonometric Fourier series, for which there is a direct connection between the differential properties of the sums of such series with the order of decreasing of the coefficients of these series.

2. Some general properties of trigonometric splines constructed by this method are investigated.

3. The trigonometric splines constructed with the help of convergence factors are investigated $\nu 1_k(r)$, $\nu 2_k(r)$ and $\nu 3_k(r)$, depending on the parameter $r$, ($r = 1, 2, ...$); It is shown that such splines are non-continuous and have continuous derivatives in order $r-1$ including.

4. Established that:

   a) trigonometric splines $ST^{(I)}(\nu 1, r, N, t)$, $ST^{(I)}(\nu 2, r, N, t)$, and $ST^{(I)}(\nu 3, r, N, t)$ coincide with each other at odd values of the parameter $r$ on grids $\Delta_N^{(I)}$, $I = 0, 1$;

   b) at even values of the parameter $r$ of trigonometric splines $ST^{(I)}(\nu 2, r, N, t)$ and coincide with each other and different from spline $ST^{(I)}(\nu 1, r, N, t)$ on grids $\Delta_N^{(I)}$, and these differences decrease rapidly with the increase of the parameter $r$;

   c) trigonometric splines $ST^{(0)}(\nu 1, r, N, t)$, $ST^{(0)}(\nu 2, r, N, t)$, and $ST^{(0)}(\nu 3, r, N, t)$ at odd values of the parameter $r$ coincide with polynomial simple splines of odd power, and splines $ST^{(1)}(\nu 2, r, N, t)$, and $ST^{(1)}(\nu 3, r, N, t)$ at even values of this parameter coincide with polynomial simple splines of even power. Establishing such a correspondence allows us to transfer to these trigonometric splines all the results of the estimates of approximation errors obtained for polynomial splines;

   d) for trigonometric splines; $ST^{(0)}(\nu 1, r, N, t)$, $ST^{(0)}(\nu 2, r, N, t)$, and $ST^{(0)}(\nu 3, r, N, t)$ at even values of the parameter $r$, and splines $ST^{(1)}(\nu 2, r, N, t)$, and $ST^{(1)}(\nu 3, r, N, t)$ at odd values of this parameter polynomial analogues are unknown; Similarly, unknown polynomial analogues to spline $ST^{(1)}(\nu 1, r, N, t)$ at all values of the parameter $r$. This implies that the class of trigonometric splines is wider than the class of polynomial splines.

5. An important positive quality of the proposed method for constructing trigonometric splines is that the complexity of constructing these splines does not depend on the specific values of the parameter $r$; and this, allows us to consider this method as an effective method for constructing polynomial simple splines of large powers.

6. Of course, trigonometric splines constructed by the proposed method require the further research.

# List of references